\documentclass[11pt]{amsart}
\usepackage{amssymb,amscd}

\newcommand{\fb}{{\mathfrak b}}
\newcommand{\fg}{{\mathfrak g}}

\newcommand{\p}{{\p}}

\def\p{\varphi}

\newtheorem{theorem}{Theorem}
\newtheorem{corollary}{Corollary}

\def\theckbibliography#1{\par\bigskip
\begin{center}
{\normalsize \bf References}
\end{center}
\par
\noindent\list
 {[\arabic{enumi}]}{\settowidth\labelwidth{[#1]}\leftmargin\labelwidth
 \advance\leftmargin\labelsep
 \usecounter{enumi}}

 \sloppy\clubpenalty4000\widowpenalty4000
 \sfcode`\.=1000\relax}

\input{tcilatex}

\begin{document}

\begin{center}
\bigskip {\large \textbf{Betti numbers of smooth Schubert varieties\ \\[0pt]
and the remarkable formula of Kostant, Macdonald\\[0pt]
Shapiro and Steinberg}}

\bigskip\textsc{Ersan Akyildiz\\[0pt]
James B.\ Carrell}
\end{center}

\noindent {\tiny \textbf{Abstract} The purpose of this note is to give a
refinement of the product formula proved in \cite{AC} for the Poincar\'e
polynomial of a smooth Schubert variety in the flag variety of an algebraic
group $G$ over ${\mathbb{C}}$. This yields a factorization of the number of
elements in a Bruhat interval $[e,w]$ in the Weyl group $W$ of $G$ provided
the Schubert variety\ associated to $w$ is smooth. This gives an elementary
necessary condition for a Schubert variety in the flag variety to be smooth.}

\bigskip

\baselineskip14pt

\section{Introduction}

Let $G$ be a semisimple linear algebraic group over ${\mathbb{C}}$, $B$ a
Borel subgroup of $G$ and $T\subset B$ a maximal torus. Let $\Phi$ denote
the root system of the pair $(G,T)$ and $\Phi^+$ the set of positive roots
determined by $B$. Let ${\alpha}_1, \dots, {\alpha}_\ell$ denote the basis
of $\Phi$ associated to $\Phi^+$, one recall the \emph{height} of ${\alpha}%
=\sum k_i{\alpha}_i \in \Phi$ is defined to be $ht({\alpha})=\sum k_i$.
Finally, let $W=N_G(T)/T$ is the Weyl group of $(G,T)$.

A remarkable formula due to Kostant \cite{KOST}, Macdonald \cite{MACD},
Shapiro and Steinberg \cite{STEIN} says that 
\begin{eqnarray}  \label{remarkform}
\sum_{w\in W} t^{2\ell(w)} = \prod_{{\alpha} \in \Phi^+}\dfrac{1-t^{2\text{ht%
}({\alpha})+2}}{1-t^{2\text{ht}({\alpha})}} = \prod_{i=1}^\ell (1+t^2+\cdots
+ t^{2m_i}),
\end{eqnarray}
where $m_1, \dots, m_\ell$ are the exponents of $G$. The left hand side of
(1) is a well known expression for the Poincar\'e polynomial\ $P(G/B,t)$ of
the flag variety $G/B$. An equivalent formulation of this identity is that
if the exponents of $(G,T)$ are ordered so that $m_\ell \ge \cdots \ge m_1$,
then the corresponding partition is dual to the partition $h_1 \ge \cdots
\ge h_k$ of $|\Phi^+|$, where $h_i$ is the number of roots of height $i$ and 
$k=h-1$ is the height of the highest root, $h$ being the Coxeter number.

In \cite{AC}, the authors gave a proof of the first half of (1) using the
fact that the cohomology algebra of $G/B$ is the coordinate algebra of the
fixed point scheme of a certain $G_a$ action on $G/B$, which is the
unipotent radical of a certain ${\mathfrak{B}}$-action as explained in
Section 2 below. Consequently, one is able to deduce a generalization of (1)
for each smooth Schubert variety\ in $G/B$. The purpose of this note is
essentially to elaborate on this aspect of (1). Along the way, we will
present an improved version of the product formula for the Poincar\'e
polynomial\ of a so called ${\mathfrak{B}}$-regular variety.

Before stating our results, let us recall the notion of a Schubert variety
in the $G/B$ setting and some of the (associated) combinatorics of $W$. For
each $w\in W$, let $wB$ denote the coset $n_wB\in G/B$, where $w=n_wT$. It
is well known that the unique $T$-fixed point in the $B$-orbit $BwB$ is $wB$%
. Furthermore, by appealing to the double coset (Bruhat) decomposition $BWB$
of $G$, $G/B$ is the union of these $B$-orbits as $w$ varies over $W$.
Furthermore, if $v$ and $w$, are distinct elements of $W$, then $vB\ne wB$,
so $W$ and $(G/B)^T$ may be identified by putting $w=wB$. Thus, $Bw$ denotes 
$BwB$. The Schubert variety\ in $G/B$ associated to $w$ is by definition the
Zariski closure $X(w)$ of $Bw$. It is well known that $\dim X(w)=\ell (w)$,
where $\ell(w)$ is the combinatorial length of $w$ with respect to the above
simple reflections $r_{{\alpha}_1}, \dots, r_{{\alpha}_\ell}$ associated to
the simple roots defined above. The Bruhat-Chevalley order $\le$ on $W$ is
the partial order on $W$ defined by putting $x\le w$ iff $X(x)\subset X(w)$.
It coincides with the standard Bruhat order defined combinatorially in terms
of the reflections. It is well known that the $Bx$ for all $x\le w$ give an
affine paving of $X(w)$. Thus the Poincar\'e polynomial\ of $X(w)$ is given
by the formula 
\begin{equation}  \label{pp}
P(X(w),t)=\sum_{x\le w} t^{2\ell(x)}.
\end{equation}

The first identity in (\ref{remarkform}) is a consequence of basic results
about smooth projective varieties admitting a regular action. In particular,
any smooth Schubert variety\ $X(w)$ in $G/B$ being such a variety, one
obtains an analogous product formula: namely, 
\begin{eqnarray}  \label{prodformsv}
P(X(w),t)= \prod_{{\alpha}\in \Phi(w)}\dfrac{1-t^{2\text{ht}({\alpha})+2}}{%
1-t^{2\text{ht}({\alpha})}},
\end{eqnarray}
where $\Phi(w)=\{{\alpha}>0\mid r_{\alpha} \le w\}$. It is well known (see 
\cite{CP} for example) that when $X(w)$ is smooth, $-\Phi(w)$ is the set of
weights of the action of $T$ on the Zariski tangent space $T_e(X(w))$ at the
identity coset.

In particular, the Poincar\'e polynomial\ of a smooth Schubert variety\
admits a factorization 
\begin{eqnarray}  \label{prodformsv2}
P(X(w),t)= \prod_{1\le i \le r} (1+t^2+\cdots +t^{2k_i})
\end{eqnarray}
for certain positive integers $k_i$. An explicit description of the $k_i$
and their multiplicities is precisely stated below in Theorem \ref{thssv}.

\section{Regular Actions and the Product Formula}

\label{RA} The product formula (\ref{prodformsv}) is a consequence of the
fact that smooth Schubert varieties\ admit a regular action. To recall this
notion, let ${\mathfrak{B}}$ denote the upper triangular Borel subgroup of $%
SL(2,{\mathbb{C}})$, and consider an algebraic action ${\mathfrak{B}}
\circlearrowright X$ on a smooth complex projective variety $X$. When the
unipotent radical ${\mathfrak{U}}$ of ${\mathfrak{B}}$ has a exactly one
fixed point $o$ on $X$, we call ${\mathfrak{B}} \circlearrowright X$ regular
and say $X$ is a ${\mathfrak{B}}$-regular variety. The following assertions
are are fundamental for our results. The first three are proved in \cite%
{CRELLE} and the last in \cite{BB}.

\medskip (1) If $X$ is ${\mathfrak{B}}$-regular, the fixed point set of a
maximal torus of ${\mathfrak{B}}$ is finite.

\medskip (2) If ${\mathfrak{T}}$ is the maximal torus on the diagonal of ${%
\mathfrak{B}}$ then $o\in X^{\mathfrak{T}}$.

\medskip (3) If ${\lambda}:{\mathbb{C}}^* \to {\mathfrak{T}}$ is the 1-psg ${%
\lambda}(s)=\text{diag}[s,s^{-1}]$, then the Bialynicki-Birula cell 
\begin{equation}  \label{bigcell}
X_o=\{x\in X\mid \lim_{s\to \infty} {\lambda}(s)\cdot x=o\}
\end{equation}
is a dense open subset of $X$. Consequently, the weights of ${\lambda}$ on $%
T_o(X)$ are all negative.

\medskip (4) $X_o$ is ${\mathfrak{T}}$-equivariantly isomorphic with the
Zariski tangent space $T_o(X)$.

Let $a_1> \dots > a_k$ denote the distinct weights of ${\lambda}$ on $T_o(X)$%
, and let $M_{a_i}$ denote the weight space corresponding to $a_i.$ Thus, 
\begin{equation}  \label{tanspace}
T_o(X)=M_{a_1} \oplus M_{a_2} \oplus \cdots \oplus M_{a_k}.
\end{equation}
As $o$ is also ${\mathfrak{B}}$-fixed, $T_o(X)$ is a $\func{Lie}({\mathfrak{B%
}})$-module. Now $\func{Lie}({\mathfrak{B}})$ is generated by 
\begin{equation*}
h={\lambda}_*(1)=%
\begin{pmatrix}
1 & 0 \\ 
0 & -1%
\end{pmatrix}
\, \,\, \text{and} \,\,\, v=%
\begin{pmatrix}
0 & 1 \\ 
0 & 0%
\end{pmatrix}%
,
\end{equation*}
and since $[h,v]=2v$, we have $hv(m)=2v(m)+vh(m)$ for all $m\in T_o(X)$.
Hence, $v(M_{a_i})\subset M_{a_i+2}.$ We will say that a regular variety is 
\emph{homogeneous} when all nonzero elements of $\ker v$ have the same
weight: equivalently, $\ker v\subset M_{a_1}$. Homogeneity has a number of
nice consequences. For example, if $i>1$, then $v$ is injective on $M_{a_i}$%
, so the weight spaces have non-increasing dimension: $\dim M_{a_j}\le \dim
M_{a_i}$ if $j>i$. The following result is proven in \cite[Theorem 3]{AC}.

\begin{theorem}
\label{multiplicity} Suppose $X$ is a homogeneous regular variety. Then:

\medskip \noindent $(a)$ $\ker v = M_{a_1}$;

\medskip \noindent $(b)$ if $i>j>0,$ $\dim M_{a_i}\le \dim M_{a_j}$;

\medskip \noindent $(c)$ the weights occur in a string: in fact, $a_i=-2i$
for each $i=1,\dots,k$.
\end{theorem}

By the theory of regular varieties, the cohomology algebra $H^*(X,{\mathbb{C}%
})$ of a ${\mathfrak{B}}$-regular variety\ is isomorphic to the graded
algebra ${\mathbb{C}}[X_o]/I$, where $I$ is a homogeneous ideal (in the
principal grading on ${\mathbb{C}}[X_o]$) with $\dim M_{a_i}$ generators of
degree $-a_i$. Let $m_i=\dim M_{a_i}$. The \emph{generalized
Kostant-Macdonald formula} (\ref{prodform}) follows from

\begin{theorem}
Suppose $X$ is a homogeneous ${\mathfrak{B}}$-regular variety. Then 
\begin{equation}  \label{prodform}
P(X, t)=\prod_{1\le i\le k} 
\begin{pmatrix}
\dfrac{1-t^{2i_i+2}}{1-t^{2i}}%
\end{pmatrix}%
^{m_i}.
\end{equation}%
.
\end{theorem}

\medskip For the details, see \cite{AC}. Putting $d_i=\dim M_{-2i} - \dim
M_{-2i-2}$ for $i=1, \dots k-1$ and $d_k=\dim M_{-2k}$ and applying Theorems
1 and 2, we get

\begin{corollary}
\label{prodformcor} Assuming $X$ is as above, 
\begin{equation}  \label{newprodform}
P(X,t)=\prod_{i=1}^k (1+t^2+\cdots +t^{2i})^{d_i}.
\end{equation}
Thus $b_2(X)=\dim M_{-2}$, and there exists a partition $b_2(X) =
\sum_{i=1}^k d_i$ such that the Euler number $\chi(X)$ is given by 
\begin{equation}  \label{euler}
\chi(X)=\prod_{i=1}^k (i+1)^{d_i}.
\end{equation}
Moreover, the partition of $\dim X$ given by 
\begin{equation}  \label{dualpart}
k= \cdots =k > (k-1) = \cdots = (k-1) > \cdots > 1 = \cdots = 1,
\end{equation}
where $i$ is repeated $d_i$ times, is dual to the partition 
\begin{equation}  \label{htpart}
\dim M_{-2}\ge \dim M_{-4}\ge \cdots \ge \dim M_{-2k}.
\end{equation}
\end{corollary}

\medskip \proof The identity (\ref{euler}) follows immediately from (\ref%
{newprodform}). The second assertion is a well known number theoretic and
combinatorial fact. See for example \cite{HUMP}. \qed

\section{$G/B$ as a homogeneous Regular Variety}

Let ${\mathfrak{g}}$, ${\mathfrak{b}}$ and ${\mathfrak{t}}$ denote the Lie
algebras of $G, B$ and $T$, and recall that $T_e(G/B)$ is naturally
isomorphic, as a ${\mathfrak{b}}$-module, to ${\mathfrak{g}}/{\mathfrak{b}}$%
. To see that $G/B$ is ${\mathfrak{B}}$-regular, let ${\alpha}_1, \dots ,{%
\alpha}_\ell$ be the simple roots for the root system $\Phi$ of $(G,T)$
determined by $B$, and select a ${\mathfrak{t}}$-weight vector $e_{{\alpha}%
_1}\in {\mathfrak{b}}$ for each simple root ${\alpha}_i$. Put 
\begin{equation*}
e=\sum_{i=1}^\ell e_{{\alpha}_i},
\end{equation*}
and let $h$ to be the unique element of ${\mathfrak{t}}$ such that ${\alpha}%
_i(h)=2$ for each $i$. Note that ${\alpha}(h)=2\text{ht}({\alpha})$ for any $%
{\alpha} \in \Phi$. Let ${\mathfrak{B}}$ denote the solvable subgroup of $B$
corresponding to the two dimensional solvable subalgebra ${\mathbb{C}}
h\oplus {\mathbb{C}} e$. Then the action ${\mathfrak{B}} \circlearrowright
G/B$, is regular with $o=e$ (under the identification of $W$ with $(G/B)^T$%
). The $h$-weight subspaces of $T_e(G/B)$ are the 
\begin{equation*}
M_{-2i}=\text{span}\{e_{-{\alpha}} \mid {\alpha}\in \Phi^+, \, {\alpha}%
(h)=-i\}
\end{equation*}
for $i=1, \dots ,k$, $k$ being the height of the longest root. Clearly $\ker
e =M_{-2}$, so ${\mathfrak{B}} \circlearrowright G/B$ is also homogeneous.
In particular, $d_i=h_i - h_{i+1}$ (i.e. the number of roots of height $i$
less the number of roots of height $i+1$).


Notice in particular that Corollary \ref{prodformcor} gives an interesting
expression for the order of the Weyl group: namely, 
\begin{equation*}
|W|=\prod_{i=1}^k (i+1)^{d_i},
\end{equation*}
where $k$ is the height of the highest root. For example, if $G=SL(n,{%
\mathbb{C}})$, then $W=S_n$, $k=n-1$ and each $d_i=1$, so 
\begin{equation*}
P(SL(n,{\mathbb{C}})/B,t)=\prod_{i=1}^{n-1} (1+t^2+\cdots +t^{2i}),
\end{equation*}
and $|W|=n!$. These facts are, of course, well known.



Since Schubert varieties in $G/B$ are $B$-stable, the smooth Schubert
varieties $X(w)$ are ${\mathfrak{B}}$-regular and homogeneous. Each $T$%
-stable line in $T_e(X(w))$ has weight ${\alpha}$ for some ${\alpha}\in
-\Phi(w),$ and thus 
\begin{equation*}
T_e(X(w))=\bigoplus_{{\alpha}\in \Phi(w)} {\mathbb{C}} e_{-{\alpha}}.
\end{equation*}
This expression immediately yields the $h$-weight space decomposition. In
fact 
\begin{equation*}
M_{-2i}=\text{span} \{e_{-{\alpha}}\mid ht({\alpha})=i, ~r_{\alpha} \le w\}.
\end{equation*}
Hence, if $k_w$ is the height of the highest root (or roots) ${\alpha}$ such
that $r_{\alpha} \le w$, then the weights of $h$ on $T_o(X(w))$ are $-2,-4,
\dots ,-2k_w$.

Applying Theorem \ref{multiplicity} in this setting gives

\begin{theorem}
\label{thssv} Let $X(w)$ be a smooth Schubert variety in $G/B$. Then $X(w)$
is a homogeneous ${\mathfrak{B}}$-regular variety\ such that the weights on $%
T_e(X(w))$ form the string of even integers, $-2 \ge -4\ge \cdots \ge -2k_w$%
. Furthermore, 
\begin{equation*}
\dim M_{-2i}=| \{e_{-{\alpha}}\mid ht({\alpha})=i, ~r_{\alpha} \le w\}|.
\end{equation*}
In particular, 
\begin{equation}  \label{PRODFORMULA}
P(X(w),t)=\prod_{i=1}^{k_w} (1+t^2+\cdots +t^{2i})^{d_i}.
\end{equation}
Consequently, 
\begin{equation}  \label{eqchi}
\chi(X(w))=|[1,w]|=\prod _{i=1}^{k_w} (i+1)^{d_i}.
\end{equation}
\end{theorem}

\medskip As far as we know, formulas (\ref{PRODFORMULA}) and (\ref{eqchi})
are new. The expression for $|[1,w]|$ provides a simple necessary condition
for the smoothness of $X(w)$ which only requires knowing the height of each $%
{\alpha} \in \Phi(w)$ and the Euler number $|[1,w]|$ of $X(w)$. Furthermore
the assertion on dual partitions of Corollary \ref{prodformcor} generalizes
to smooth Schubert varieties\ a well known fact about the heights of
positive roots; namely that the partition of $\dim G/B=|\Phi^+|$ given by
the exponents of $G$ is dual to the partition $h_1>h_2> \cdots >h_k$, where $%
h_i=|\{{\alpha}>0 \mid ht({\alpha})=i\}|$.

\begin{corollary}
Suppose $X(w)$ is a smooth Schubert variety, and let $h_{w,i}=|\{{\alpha}>0
\mid ht({\alpha})=i, ~ r_{\alpha} \le w\}|$. Then the partition $h_{w,1}\ge
h_{w,2} \ge \cdots \ge h_{w,k}$ of $\ell(w)$, is dual to the partition to
the partition (\ref{dualpart}) where each $i$, $1\le i \le k$, is repeated $%
d_i$ times.
\end{corollary}

\section{The formula in the rationally smooth\ case}

Suppose $X(w)$ is a Schubert variety\ for which $P(X(w),t)$ has the form (%
\ref{newprodform}) but doesn't necessarily arise from a product formula as
in (\ref{prodformsv}). Then clearly $P(X(w),t)$ is palindromic, so $X(w)$
has to be rationally smooth, by the main result of \cite{CP}. On the other
hand, one can ask if there also exist singular Schubert varieties\ for which
(\ref{prodformsv}) holds. The converse is answered by a result of Sara
Billey \cite{BILLEY}: If $X(w)$ is a rationally smooth\ Schubert variety\ in 
$G/B$, there exist not necessarily distinct positive integers $j_1, j_2,
\dots, j_k$ such that 
\begin{equation}  \label{pprssv}
P(X(w),t)=\prod_{1\le i \le k}(1+t^2+\cdots +t^{2j_i}).
\end{equation}

Perhaps as one should expect, the product formulas (\ref{pprssv}) and (\ref%
{prodformsv}) needn't coincide, except in the smooth case. To see this,
consider the following example.

\medskip \noindent \textbf{Example 1.} Let $\Phi(B_2)$ be the root system of
type $B_2$ with two simple roots ${\alpha}$ and ${\beta}$ where ${\alpha}$
is long and ${\beta}$ is short, and positive roots ${\alpha}, \, {\beta}, \, 
{\alpha}+{\beta}$ and ${\alpha}+2{\beta}$. All Schubert varieties\ in $B_2/B$
are rationally smooth, and the unique singular Schubert variety\ is $X(w)$,
where $w=r_{\beta} r_{{\alpha}+2{\beta}} r_{\beta}.$ The fact that $X(w)$ is
singular follows easily since if it were smooth, the heights 1 and 3 of ${%
\alpha},~{\beta},~{\alpha}+2{\beta}$ would have to be consecutive integers,
which obviously they aren't. The right hand side of product formula (\ref%
{prodformsv}) is 
\begin{equation*}
\frac{(1+t)(1+t)(1-t^4)}{(1-t^3)},
\end{equation*}
which isn't even a polynomial. In this case, the factorization (\ref{pprssv}%
) of $P(X(w),t)$ is 
\begin{equation*}
P(X(w),t)=(1+t)(1+t+t^2).
\end{equation*}

\medskip

In types $ADE$, (\ref{pprssv}) follows from (\ref{prodformsv}) and
Peterson's $ADE$ Theorem \cite{CK} which says that if $G$ is simply laced,
every rationally smooth\ Schubert variety\ in $G/B$ is smooth. The proof of (%
\ref{pprssv}) for types $B$ and $C$ is given in \cite{BILLEY}. The case of $%
G_2$ was checked by hand, while the $F_4$ case was verified by computer.

The above comments lead naturally to the question of which rationally
smooth\ Schubert varieties\ are smooth. The answer is quite easy to state.
Let $TE(X(w))$ be the span of the tangent lines of the $T$-stable curves in $%
X(w)$ containing $e$: that is,%
\begin{equation*}
TE(X(w))=\bigoplus_{{\alpha}\in \Phi(w)} {\mathbb{C}} e_{-{\alpha}}.
\end{equation*}

The following answer is given in (see \cite{C2}):

\begin{theorem}
\label{thm4} Suppose $G$ doesn't contain $G_2$-factors and let $X(w)$ be a
rationally smooth\ Schubert variety\ in $G/B$. Then $X(w)$ is smooth if and
only if $TE(X(w))$ is a $B$-submodule of $T_e(X(w))$.
\end{theorem}

In fact, Theorem \ref{thm4} fails in type $G_2$: there exists a singular
rationally smooth\ Schubert variety\ $X(w)$ in $G_2/B$ such that $TE(X(w))$
is a $B$-submodule of $T_e(X(w))$. The Poincar\'e polynomial\ of this
Schubert variety\ is also given by the product formula (\ref{prodformsv}).
To see this explicitly, let ${\alpha}$ and ${\beta}$ denote respectively the
long and short simple roots for $G_2$ corresponding to $B$, and let $r=r_{%
\alpha}$ and $s=r_{\beta}$ be the corresponding reflections. Let $w=srsrs$.
Now $\ell(w)=5$ and it is not hard to see that 
\begin{equation*}
\Phi(w)=\{{\alpha},{\beta},{\alpha}+{\beta},{\alpha}+2{\beta}, {\alpha}+3{%
\beta}\}.
\end{equation*}
It follows that $TE(X(w))$ is indeed a $B$-submodule of $T_e(X(w))$.
However, it is well known that $X(w)$ is singular: for example; for example,
see \cite{KUM}. Curiously, (\ref{prodformsv}) holds for $X(w)$ too. Indeed,
the heights in $\Phi(w)$ are 1, 2, 3, and 4, while $d_1=1$, $d_2=d_3=0$ and $%
d_4=1$. Thus the right hand side of (\ref{prodformsv}) is 
\begin{equation*}
(1+t^2)(1+t^2 +t^4 +t^6 +t^8),
\end{equation*}
which is indeed the Poincar\'e polynomial\ of $X(w)$. Thus, there exist
singular examples where the product formula makes sense and coincides with
Billey's factorization.

It would be interesting to know whether or not Theorem \ref{thm4} holds
under the weaker condition that $TE(X(w))$ is only a ${\mathfrak{B}}$%
-submodule. We conclude with an example illustrating Theorem \ref{thssv}.

\medskip \noindent \textbf{Example 2.} Let ${\alpha}_1, \dots ,{\alpha}_4$
denote the simple roots for the root system of $D_4$ labelled according to
the usual labelling of the Dynkin diagram \cite{HUMP}, and denote the
corresponding reflections by $1,\, 2,\, 3\, ,\, 4$. Consider the element $%
w=2142132$ of $W(D_4)$ of length 7. By Goresky's extremely useful tables 
\cite{G}, 
\begin{equation*}
P(X(w),t)=1+4t^2+9t^4+13t^6+13t^8+9t^{10}+4t^{12}+t^{14},
\end{equation*}
so $X(w)$ is smooth since its Poincar\'e polynomial\ is palindromic. One
easily checks that 
\begin{equation*}
\Phi(w) =\{1,\, 2, \,3, \,4, \, 212, \, 232, \, 242\},
\end{equation*}
so the heights are 1 and 2 with multiplicities 4 and 3 respectively. Thus, $%
d_1=1$ and $d_2=3$. Hence, Theorem \ref{thssv} gives the expression 
\begin{equation*}
P(X(w),t)=(1+t^2)(1+t^2+t^4)^3,
\end{equation*}
which agrees with the above expression.

\newpage {\footnotesize 
\begin{theckbibliography}{00}
\bibitem{AC} E. Akyildiz and J.B. Carrell: {\em A generalization of the Kostant--Macdonald
identity},  Proc. Nat. Acad. Sci. U.S.A. {\bf 86} (1989), 3934--3937.
\bibitem{BB} A. Bialynicki-Birula: {\em Some theorems on actions of algebraic groups}. Ann. of
Math. (2) {\bf 98} (1973), 480--497.
\bibitem{BILLEY} S.\ Billey: {\em Pattern avoidance and rational smoothness of Schubert
varieties}, Adv. Math. {\bf 139} (1998), no. 1, 141--156.
\bibitem{CP} J. B.\ Carrell: {\em The Bruhat Graph of a Coxeter Group, a
Conjecture of Deodhar, and Rational Smoothness of Schubert Varieties,}
Proc.\ Symp.\ in Pure Math.\ {\bf 56}, No.\ 2, (1994), Part 1, 53-61.
\bibitem{CRELLE} J. B.\ Carrell: {\em Deformation of the nilpotent zero scheme and the
intersection ring of invariant subvarieties}. J. Reine Angew. Math. {\bf 460} (1995), 37--54.
\bibitem{C} J. B.\ Carrell: {\em  $B$-submodules of $\fg/\fb$ and smooth Schubert varieties
in $G/B$}, 	arXiv:1007.4365v1.
\bibitem{CK}  J. B.\ Carrell and J. \ Kuttler: {\em Singular points of $T$-varieties in $G/P$
and the Peterson map}, Invent. Math. {\bf 151}  (2003), 353--379.
\bibitem{G} M.\ Goresky: http://www.math.ias.edu/~goresky/D4.txt.
\bibitem{HUMP} J.\ Humphreys: {\em Reflection groups and Coxeter groups}.
 Cambridge Studies in
Advanced Mathematics {\bf 29}. Cambridge University Press, Cambridge, 1990.
\bibitem{KOST} B.\ Kostant: {\em The principal three-dimensional subgroup 
and the Betti numbers of a complex simple Lie group}, Amer. J. Math. 
{\bf 81} (1959), 973--1032. 
\bibitem{KUM} S.\ Kumar: {\em  Nil Hecke ring and singularity of Schubert
varieties}, Inventiones Math. {\bf 123} (1996), 471--506.
\bibitem{MACD} I. G.\ Macdonald: {\em The Poincar\'e series of a Coxeter group},  Math. Ann.  
{\bf 199} (1972), 161--174. 
\bibitem{STEIN} R.\ Steinberg: {\em Finite reflection groups},
Trans. Amer. Math. Soc. {\bf 91} (1959), 493--504. 
\end{theckbibliography}
}

{\footnotesize \medskip \noindent \noindent {\tiny Ersan Akyildiz\newline
Institute of Applied Mathematics \newline
Middle East Technical University\newline
Ankara 06531, Turkey\newline
ersan$@$metu.edu.t{r}} }

{\footnotesize \medskip \noindent {\tiny James B.\ Carrell \newline
Department of Mathematics \newline
University of British Columbia \newline
Vancouver, Canada V6T 1Z2 \newline
carrell$@$math.ubc.c{a}} }

\end{document}